\documentclass[dvipdfmx]{amsart}
\usepackage[cmtip,arrow]{xy}
\usepackage{
amsxtra,amsmath,amsthm,amssymb,amscd,ascmac,amsfonts,
multicol,tabularx,latexsym,mathrsfs,
pb-diagram,pb-xy
}
\usepackage{mathtools} 
\usepackage{amsfonts}
\theoremstyle{definition}
\usepackage{tikz}
\usetikzlibrary{cd, decorations.pathmorphing}
\newcommand{\longhookrightarrow}{\begin{tikzcd}[cramped,sep=scriptsize,ampersand replacement=\&]{}\arrow[r, hook]\&{}\end{tikzcd}} 
\newcommand{\longtwoheadrightarrow}{\begin{tikzcd}[cramped,sep=scriptsize,ampersand replacement=\&]{}\arrow[r, two heads]\&{}\end{tikzcd}} 

\newcommand{\Q}{\mathbb Q}
\newcommand{\Z}{\mathbb Z}

\newcommand{\R}{\mathbb R}
\newcommand{\C}{\mathbb{C}}
\newcommand{\Gal}{\mathrm{Gal}}
\newtheorem{thm}{Theorem}

\newtheorem{thma}{Theorem\,A\!\!}
\newtheorem{thmb}{Theorem\,B\!\!}
\newtheorem{thmc}{Theorem\,C\!\!}
\newtheorem{thmd}{Theorem\,D\!\!}
\newtheorem{thme}{Theorem\,E\!\!}
\newtheorem{lem}{Lemma}

\newtheorem{prop}{Proposition}
\newtheorem{rem}{Remark}

\newcommand{\F}{\mathbb F}
\newcommand{\Fb}{\overline{\mathbb{F}}}

\newcommand{\ab}{\mathrm{ab}}

\newdimen\minCDarrowwidth
\date{}
\dedicatory{Dedicated to the memory of Tamao Ozaki}
\begin{document}
\title[]
{A number field analogue of the Grothendieck conjecture for
curves over finite fields}
\author{Manabu\ OZAKI}
\maketitle
\section{Introduction}
Analogy between number fields and
1-dimensional function fields (or algebraic curves) over finite fields has led us to a deep insight into these two arithmetic objects.
For example, the ``Main Conjecture" of Iwasawa theory
of cyclotomic $\Z_p$-extensions
can be regard as an analogy of Weil's
theorem on the relationship
among the congruent zeta function and the Frobenius action on the Tate module associated to a curve over a finite field.
Here, we regard that the cyclotomic $\Z_p$-extension, or adjoining all the $p$-power-th roots of unity, is analogous
to the constant field extension $\Fb K/K$
of a function field $K$ over a finite
field $\F$ to an algebraic closure $\Fb$.

However the extension $\Fb K/K$
is in fact given by
ajoining {\it all} the roots of unity in $\overline{\F}$.
Therefore it is natural to ask what
happens if we consider the {\it maximal cyclotomic extension} of a number field $k$, namely, the extension $k(\mu_\infty)/k$
given by
adjoining all the roots of unity $\mu_\infty$,
instead of the cyclotomic $\Z_p$-extension.

In the present paper,
we choose the maximal cyclotomic extension
as the analogous object of the constant field extension $\Fb K/K$.
Then we will give
a number field analogue of the Grothendieck conjecture
for curves over finite fields,
proved by Tamagawa \cite{Tam} and 
Mochizuki \cite{Moch07}.

Let $C$ be a non-singular geometrically connected curve over a
field $F$.
Denote by $F(C)$ the function field of $C$ over $F$,
and by $\overline{F}^{\mathrm{sep}}$ a separable closure of $F$.
We put $\overline{F}^{\mathrm{sep}}(C):=F(C)\overline{F}^{\mathrm{sep}}$,
and
define $L(C)$ to be the maximal extension field (in a fixed algebraic closure) of 
$\overline{F}^{\mathrm{sep}}(C)$
unramified at 
$C\times_F\overline{F}^{\mathrm{sep}}$
(Note that $L(C)/\overline{F}^{\mathrm{sep}}(C)$ is the maximal unramified extension if $C$ is projective).

Then we get the following fundamental exact sequence:
\begin{equation}\label{fe}
1\longrightarrow\Gal(L(C)/\overline{F}^{\mathrm{sep}}(C))
\longrightarrow\Gal(L(C)/F(C))
\overset{P_C}{\longrightarrow}G_F
\longrightarrow 1,
\end{equation}
where 
$G_F:=\Gal(\overline{F}^{\mathrm{sep}}/F)\simeq\Gal(\overline{F}^{\mathrm{sep}}(C)/F(C))$
is the absolute Galois group of $F$.

Assume that $C$ is
{\it hyperbolic}, namely,
$2-2g_C-n_C<0$ holds,
where $g_C$ is the genus of $C$ and 
$n_C:=\#(C^*(\overline{F})-C(\overline{F}))$, 
$C^*$ being the compactification of $C$.
This assumption is equivalent to that $\Gal(L(C)/\overline{F}^{\mathrm{sep}}(C))$ is non-abelian
if $\mathrm{char}(F)=0$.
The Grothendieck conjecture asserts that
the pro-finite group homomorphism $P_C:\Gal(L(C)/F(C))\longtwoheadrightarrow G_F$ has
much information to reconstruct
the curve $C$ itself.
This conjecture has been established by A.Tamagawa
(case $n_C>0$)
and S.Mochizuki (case $n_C=0$) in the case where
$F$ is finitely generated over $\Q$ or
finite.
More presisely:
\begin{thma}[Tamagawa\cite{Tam},  Mochizuki\cite{Moch03}]
Let $F$ be a field finitely generated over $\Q$, and $C_i$ a hyperbolic curve over $F$
($i=1,2$). Put
$G_i:=\Gal(L(C_i)/F(C_i))$ ($i=1,2$).
Then there is the natural bijection
\[
\mathrm{Isom}_F(C_1,C_2)
\simeq
\mathrm{Isom}_{G_F}(G_1,G_2)
/\mathrm{Inn}(\Gal(L(C_2)/\overline{F}(C_2)),
\]
where the left hand side is the set
of the $F$-isomorphisms of curves from $C_1$ to
$C_2$ and the right hand side is
the quotient of the set of the pro-finite group isomorphisms
from
$G_1$ to $G_2$
which are compatible with $P_{C_1}$ and
$P_{C_2}$
by the right action of the inner automorphism group of $\Gal(L(C_2)/\overline{F}(C_2))$. 
\end{thma}
In the case where $F$ is finite, the following
``absolute version" of the Grothendieck
conjecture holds:
\begin{thmb}[Tamagawa\cite{Tam}, Mochizuki\cite{Moch07}]
Let $F_i$ be a finite field, and $C_i$
a hyperbolic curve over $F_i$ ($i=1,2$).
Put $G_i:=\Gal(L(C_i)/F_i(C_i))$
($i=1,2$).
Then there is the natural bijection
\[
\mathrm{Isom}(C_1,C_2)
\simeq
\mathrm{Isom}(G_1,G_2)
/\mathrm{Inn}(G_2),
\]
where the left hand side is the set
of the isomorphisms as schemes from $C_1$ to
$C_2$ and the right hand side is
the quotient of the set of the pro-finite group isomorphisms
from
$G_1$ to $G_2$
by the right action of the inner automorphism group of $G_2$. 
\end{thmb}
From the viewpoint that 
the maximal cyclotomic extension $k(\mu_\infty)/k$ of a number field $k$ is an analogous object to the constant field extension $\Fb(C)/\F(C)$ for the function field of a curve $C$ over a finite field $\F$,
we find that the number field analogue of the fundamental exact sequence \eqref{fe} in the case where $C$ is projective
is
\[
1\longrightarrow\Gal(L(\tilde{k})/\tilde{k})
\longrightarrow\Gal(L(\tilde{k})/k)
\longrightarrow\Gal(\tilde{k}/k)
\longrightarrow 1,
\]
where $L(\tilde{k})$ is the maximal unramified exttension over $\tilde{k}:=k(\mu_\infty)$.

Under this situation, we will give the following theorem analogous to Theorem B:
\begin{thm}
Let $k_i$ be a number field of finite degree,
and $L_i$ the maximal unramified extension field over $k_i(\mu_\infty)$ ($i=1,2$).
For any isomorphism
\[
\varphi:\Gal(L_1/k_1)\simeq\Gal(L_2/k_2)
\]
of pro-finite groups, 
there exists the unique field isomorphism
\[
\tau:L_1\simeq L_2
\]
such that
$\tau(k_1)=k_2$ and
\[
\varphi(x)=\tau x\tau^{-1}
\]
holds for every $x\in\Gal(L_1/k_1)$.
In other words, we have the natural bijection
\[
\mathrm{Isom}(L_1/k_1,L_2/k_2)
\simeq
\mathrm{Isom}(\Gal(L_1/k_1),\Gal(L_2/k_2)),
\]
where the left hand side is the set of
all the field isomorphisms $\tau:L_1\simeq L_2$ with $\tau(k_1)=k_2$, and the right hand side is the set of all the isomorphism
$\Gal(L_1/k_1)\simeq\Gal(L_2/k_2)$
of pro-finite groups.
\end{thm}
\begin{rem}
The above theorem is an affirmative
answer to \cite[Conjecture(12.5.3)]{NSW}
in the case where the ramified prime sets
are empty.
\end{rem}
\section{Arithmetically equivalence}
For any number field $F$,
we put $\tilde{F}:=F(\mu_\infty)$,
and denote by $L(F)$ the maximal unramified extension field over $F$. 

In this section,
we will show the following proposition,
which plays a crucial role in the proof of Theorem 1.
\begin{prop}\label{prop1}
Let $k_i$ be a number field of finite degree
($i=1,2$).
Assume that there exists an isomorphism
\[
\varphi:\Gal(L_1/k_1)\simeq\Gal(L_2/k_2)
\]
of pro-finite groups, where $L_i:=L(\tilde{k_i})\ (i=1,2)$.

Then for any finite subextension $M_1/k_1$ of $L_1/k_1$, $M_1$ and $M_2:=L_2^{\varphi(\Gal(L_1/M_1))}$
are arithmetically equivalent:
$M_1\approx M_2$.
Namely, they have the same
Dedekind zeta function: $\zeta_{M_1}(s)=\zeta_{M_2}(s)$.
\end{prop}
In Proposition \ref{prop1},
we note that $L_i=L(\tilde{M_i})$ and
$\varphi$ induces $\Gal(L_1/M_1)\simeq
\Gal(L_2/M_2)$,
hence it is enough to show the proposition in the case where
$M_1=k_1$.

Let $K_i/k_i$ be the maximal abelian subextension of $L_i/k_i$, 
and denote by $L^\mathrm{ab}(K_i)$
the maximal unramified abelian extension field over $K_i$.
Since $L^\mathrm{ab}(K_i)\subseteq L_i$
and
\begin{align*}
\Gal(L^\mathrm{ab}&(K_i)/k_i)\\
&\simeq
\Gal(L_i/k_i)/((\Gal(L_i/k_i),\Gal(L_i/k_i)),(\Gal(L_i/k_i),\Gal(L_i/k_i))),
\end{align*}
Proposition \ref{prop1} follows from the following:
\begin{prop}\label{prop2}
Assume that there exists an isomorphism
\[
\psi:\Gal(L^{\mathrm{ab}}(K_1)/k_1)\simeq
\Gal(L^\mathrm{ab}(K_2)/k_2)
\]
of pro-finite groups.
Then we have $k_1\approx k_2$.
\end{prop}
Our proof of the above proposition is based on the following fact:
\begin{thmc}[Stuart-Perlis\cite{S-P}]
Let $g_F(l)$ be the number of primes of $F$ lying over $l$ for any number field $F$
of finite degree and prime number $l$.
Assume that $g_{k_1}(l)=g_{k_2}(l)$ holds
for number fields $k_1$ and $k_2$ of finite degree and all but finitely many prime numbers $l$.
Then we have $k_1\approx k_2$.
\end{thmc}
In what follows, we will show that
$g_{k_i}(l)$ is encoded in the group structure of $\Gal(L^\mathrm{ab}(K_i)/k_i)$.
For simplicity, we write $K$ and $k$
for $K_i$ and $k_i$, respectively.

We first analyze $\Gal(K/k)$
and its decomposition and inertia subgroups.
We denote by $D_{\frak{l}}(M/F)$ and 
$I_\frak{l}(M/F)$ the decomposition 
and the inertia subgroup, respectively,
of $\Gal(M/F)$ for a prime $\frak{l}$ of $k$,
$M/F$ being any abelian extension of number fields.
\begin{lem}\label{lem-DI}
For any non-archimedean prime $\frak{l}$ of $k$, we have
\[
D_\frak{l}(K/k)\simeq\hat{\Z}\times\Z_l
\times\Z/d_\frak{l},\ \ \ 
I_{\frak{l}}(K/k)\simeq\Z_l\times
\Z/d_\frak{l},
\]
where $l$ is the rational prime below $\frak{l}$ and $d_\frak{l}$ is 
a certain divisor of $l-1$ (if $l\ne 2$) 
or 2 (if $l=2$).
Furthermore, we have $d_\frak{l}=l-1$ (if $l\ne 2$) or $d_\frak{l}=2$ (if $l=2$))
if $\frak{l}$ is unramified in $k/\Q$.
\end{lem}

{\bf Proof.}\ \ \ 
For any rational prime $l$, it is easy to see that
\[
D_l(\tilde{\Q}/\Q)\simeq\hat{\Z}\times\Z_l
\times\Z/b_l,\ \ \ 
I_l(\tilde{\Q}/\Q)\simeq\Z_l\times
\Z/b_l,
\]
where $b_l=l-1$ (if $l\ne 2$) or $b_l=2$
(if $l=2$).

The restriction $\Gal(\tilde{k}/k)
\longhookrightarrow\Gal(\tilde{\Q}/\Q)$
induce injections 
$D_\frak{l}(\tilde{k}/k)
\longhookrightarrow
D_l(\tilde{\Q}/\Q)$
and $I_\frak{l}(\tilde{k}/k)
\longhookrightarrow I_l(\tilde{\Q}/\Q)$
with finite cokernels.
Hence we see that
\begin{equation}\label{di}
D_\frak{l}(\tilde{k}/k)\simeq\hat{\Z}\times\Z_l
\times\Z/d_\frak{l},\ \ \ 
I_{\frak{l}}(\tilde{k}/k)\simeq\Z_l\times
\Z/d_\frak{l},
\end{equation}
for a certain divisor $d_l$ of $l-1$ or $2$,
which equals $l-1$ or $2$ if $l$ is unramified in $k/\Q$.
Because $K/\tilde{k}$ is unramified and
every prime totally splits in it,
we see that the restriction
$\Gal(K/k)\longtwoheadrightarrow\Gal(\tilde{k}/k)$
induces isomorphisms
\[
D_\frak{l}(K/k)\simeq D_\frak{l}(\tilde{k}/k),\ \ 
I_\frak{l}(K/k)
\simeq I_\frak{l}(\tilde{k}/k).
\]
Therefore the assertion of the lemma follows
from
\eqref{di}.
\hfill$\Box$
\begin{lem}\label{lem-sep}
$I_{\frak{l_1}}(K/k)\cap
I_{\frak{l}_2}(K/k)=1$ for any non-archimedean primes $\frak{l}_1$ and $\frak{l}_2$ of $k$ with $\frak{l}_1\ne\frak{l}_2$. 
\end{lem}
{\bf Proof.}\ \ \ 
We first note that $K$ is the genus class field of $\tilde{k}/k$, namely,
the maximal unramified abelian extension filed of $\tilde{k}$ which is abelian
over $k$.
Then, by using \cite[Proposition 2]{Fur},
we see that the global resiprocity map induces
the isomorphism
\begin{equation}\label{rec}
\rho:\mathcal{C}_k:=J_k\left/\overline{k^\times
\prod_{\frak{p}\mbox{\tiny : primes of $k$}}U^0_{\frak{p}}}\right.
\simeq\Gal(K/k),
\end{equation}
where $J_k$ is the idele group of $k$,
\[
U_{\frak{p}}^0:=
\begin{cases}
\ker(N_{k_\frak{p}/\Q_p}
\!:U_{\frak{p}}\longrightarrow\Z_p^\times)\ \mbox{if $\frak{p}$ is a non-archimedean prime}, 
\\
\R^\times_{>0}\ \mbox{if $\frak{p}$ is a real arichimedean prime},
\\
\C^\times\ \mbox{if $\frak{p}$ is a complex
archimedean prime},
\end{cases}
\]
$U_{\frak{p}}$ and $p$
being the local unit group of $k_{\frak{p}}$
and the prime number below $\frak{p}$, respectively, 
and the ``bar" means topological closure.

Assume that $I_{\frak{l}_1}(K/k)
\cap I_{\frak{l}_2}(K/k)\ne 1$.
Then it follows from \eqref{rec} that
there exists
$u_{\frak{l}_i}\in U_{\frak{l}_i}\ (i=1,2)$
such that 
$[(u_{\frak{l}_1})]=[(u_{\frak{l}_2})]\ne 1$,
where $[(u_{\frak{l}_i})]\in\mathcal{C}_k$
stands for the image of $u_{\frak{l}_i}\in
U_{\frak{l}_i}$
under the composite of the natural maps
\[
U_{\frak{l}_i}
\longhookrightarrow J_k
\longtwoheadrightarrow \mathcal{C}_k,
\ \ 
u_{\frak{l}_i}\mapsto(u_{\frak{l}_i})
\mapsto[(u_{\frak{l}_i})].
\]
Suppose that $\frak{l}_1\ne\frak{l}_2$.
Then we have
\begin{equation}\label{closure}
1\ne(u_{\frak{l}_1})(u_{\frak{l}_2})^{-1}
\in
\overline{k^\times
\prod_{\frak{p}\mbox{\tiny : primes of $k$}}U^0_{\frak{p}}}.
\end{equation}
It follows from \cite[ Théorème 1]{Ch} that for any given integer $m\ge 1$, there exists
a finite set $T_m$ of degree one primes of $k$
such that $\frak{l}_1,\frak{l}_2\not\in T_m$ and
\[
E_k\ni\varepsilon\equiv 1\pmod{\prod_{\frak{p}\in T_m}\frak{p}}
\Longrightarrow \varepsilon\in E_k^m
\]
holds, $E_k$ being the grobal unit group of $k$.
Because $U_\frak{p}^0=1$ for $\frak{p}\in T_m$, we find that
there exists an open neighberhood $U_k\supseteq H_m\ni(u_{\frak{l}_1})(u_{\frak{l}_2})^{-1}$,
$U_k\subseteq J_k$ being the unit idele group of $k$,
such that 
\[
H_m\cap k^\times
\prod_{\frak{p}\mbox{\tiny : primes of $k$}}U^0_{\frak{p}}\subseteq E_k^m\prod_{\frak{p}\mbox{\tiny : primes of $k$}}U^0_{\frak{p}}.
\]
Hence, by \eqref{closure}, we see that for each $n\ge 1$,
there exists 
$\varepsilon_n\in E_k$ and 
$w_{n,\frak{l}_i}\in U_{\frak{l}_i}^0$
such that
\[
u_{\frak{l}_1}\equiv
\varepsilon_n
^{\#\left(\mathcal{O}_{\frak{l}_1}/\frak{l}_1^n\mathcal{O}_{\frak{l}_1}\right)^\times}
w_{n,\frak{l}_1}
\equiv w_{n,\frak{l}_1}\pmod{\frak{l}_1^n\mathcal{O}_{\frak{l}_1}},
\]
$\mathcal{O}_{\frak{l}_1}$ being the integer ring of $k_{\frak{l}_1}$,
which implies $u_{\frak{l}_1}\in  U^0_{\frak{l_1}}$ since $U^0_{\frak{l_1}}$
is closed in $U_{\frak{l}_1}$.
This contradicts to 
$[(u_{\frak{l}_1})]\ne 1$.
Thus we conclude that $I_{\frak{l}_1}(K/k)
\cap I_{\frak{l}_2}(K/k)\ne 1$
implies $\frak{l_1}=\frak{l_2}$.
\hfill$\Box$

\

\begin{lem}\label{lem-un}
Let $M$ be an intermediate field of $K/k$ such that
for each finite subextension $F/k$ of $M/k$, there exist
infinitely many degree one primes $\mathfrak{L}$ of $F$ such that
$\mu_l\subseteq M$, $l$ being the rational prime below $\mathfrak{L}$.
Then we have 
\[
\varprojlim_{k\subseteq F\subseteq M,\,[F:k]<\infty} E_F=0,
\]
where $E_F$ denotes the global unit group of $F$ and
the projective limit is taken with respect to the norm maps.
\end{lem}
{\bf Proof.}\ \ \ 
It is enough to show that 
\[
\bigcap_{F\subseteq N\subseteq M,\ [N:F]<\infty} N_{N/F}(E_N)=1
\]
for each finite subextension $F/k$ of $M/k$.
Let $\frak{L}$ be a degree one prime of $F$ such that
$\mu_l\subseteq M$, $l$ being the rational prime below $\mathfrak{L}$,
and $\frak{L}$ is unramified in $F/\Q$.
Then $F(\mu_l)\subseteq M$ and 
$N_{F(\mu_l)/F}(\varepsilon)\equiv 1\pmod{\frak{L}}$ holds
for every $\varepsilon\in E_{F(\mu_l)}$.
Since $F(\mu_l)\subseteq M$, we find that
$\eta\equiv 1\pmod{\frak{L}}$ holds for any
$\eta\in\bigcap_{F\subseteq N\subseteq M,\ [N:F]<\infty}
N_{N/F}(E_N)$.
Because there exists infinitely many such primes $\frak{L}$,
we conclude that $\eta=1$ must hold.
\hfill$\Box$

\

\begin{lem}\label{lem-Y}
Let $p$ be an odd prime number and $\Delta=\langle\delta\rangle\subseteq\Gal(K/k)$ a subgroup of order $p$.
Put 
\[
Y_\Delta^{(p)}:=\Gal(L_p^{\mathrm{ab}}(K)/K^{\Delta})
/(\Gal(L_p^{\mathrm{ab}}(K)/K),\overline{\delta}),
\]
where 
$L_p^{\mathrm{ab}}(K)/K$ is the maximal $p$-subextension 
of $L^{\mathrm{ab}}(K)/K$, and
$\overline{\delta}$ is a lift of $\delta$ to $\Gal(L_p(K)/K^\Delta)$.
Then
\[
\mathrm{Tor}(Y_{\Delta}^{(p)})
\simeq
\begin{cases}
0\ \ \mbox{if $\Delta\not\subseteq I_{\frak{l}}(K/k)$ for any prime $\frak{l}$
of $k$},\\
\F_p[[\Gal(K/k)/D_{\frak{l}}(K/k)]]
\ \ 
\mbox{if $\Delta\subseteq I_{\frak{l}}(K/k)$
for some prime $\frak{l}$ of $k$},
\end{cases}
\]
as $\F_p[[\Gal(K/k)]]$-modules,
where $\mathrm{Tor}(Y_\Delta^{(p)})$ means the torsion part
of the pro-$p$ abelian group $Y_\Delta^{(p)}$.
\end{lem}

{\bf Proof.}\ \ \ 
Let $M
:=L^{\mathrm{ab}}_p(K)^{(\Gal(L^{\mathrm{ab}}_p(K)/K),\overline{\delta})}$.
Then $M/K^\Delta$ is the maximal abelian $p$-subextension of $L^{\mathrm{ab}}_p(K)/K^\Delta$,
and $\Gal(K/k)$ acts on
$Y_\Delta^{(p)}=\Gal(M/K^\Delta)$ via inner automorphisms
of $\Gal(M/k)$.

Assume that $\Delta\not\subseteq I_{\frak{l}}(K/k)$ for any prime $\frak{l}$ of $k$.
Then $M$ is the maximal unramified abelian $p$-extension field over $K^{\Delta}$ since
$K/K^\Delta$ is unramified.

Now we employ the following theorem:
\begin{thmd}[Uchida\cite{Uch82}]
For any integer $m\ge 1$ and prime number $l\equiv 1\pmod{m}$, 
denote by $\Q(l,m)$ the subfield
of the $l$-th cyclotomic field $\Q(\mu_l)$
such that $[\Q(\mu_l):\Q(l,m)]=m$.
We define $\Q^{(m)}$ to be the composite
field of all the $\Q(l,m)$ for the primes
$l\equiv 1\pmod{m}$.

Let $F$ a number field with $\Q^{(m)}\subseteq F$ for
some $m\ge 1$.
Furthermore we assume that
$F$ contains a subfield $F_0$ of finite
degree over $\Q$ such that $F$ is a subfield of the maximal nilpotent extension
of $F_0$.
Then the Galois group
the maximal unramified pro-solvable extension
over $F$ is a free pro-solvable group
of countably infinite rank.
\end{thmd}

\

It follows from Theorem D that
$Y_{\Delta}^{(p)}$ is a free pro-$p$ abelian group
since $\Q^{(p)}\subseteq K^{\Delta}$ and
$K^\Delta/k$ is abelian.
Hence $\mathrm{Tor}(Y_\Delta^{(p)})=0$ in this case.

Assume that $\Delta\subseteq I_\frak{l}(K/k)$
for a certain non-archimedean prime $\frak{l}$
of $k$.
Then such a prime $\frak{l}$ is unique by 
Lemma \ref{lem-sep}, and 
the prime number $l$ below $\frak{l}$
satisfies $l\equiv 1\pmod{p}$ by Lemma
\ref{lem-DI}, especially, $K/K^\Delta$
is tamely ramified.
We get the exact sequence
\[
1\longrightarrow\langle I_\frak{L}(M/K^\Delta)\,
|\ \frak{L}|\frak{l}\,\rangle
\longrightarrow
\Gal(M/K^{\Delta})
\longrightarrow
\Gal(L_p^{\mathrm{ab}}(K^\Delta)/K^{\Delta})
\longrightarrow 1
\]
of abelian pro-$p$-groups,
where $\frak{L}$ runs over all the primes of $K^\Delta$ lying over $\frak{l}$.
Here we note that exactly all the primes lying over $\frak{l}$ ramify in $M/K^\Delta$,
and that $L^\mathrm{ab}(K^\Delta)K\subseteq M$
because $M/K$ is the maximal unramified abelian
$p$-extension which is abelian over $K^\Delta$.

It follows from Theorem D that
$\Gal(L^{\mathrm{ab}}(K^\Delta)/K^{\Delta})$
is free pro-$p$ abelian group since
$\Q^{(p)}\subseteq K^{\Delta}$, hence the above exact sequence splits.
Thetrefore we obtain
\begin{equation}\label{tor}
\mathrm{Tor}(Y_\Delta^{(p)})=\langle I_\frak{L}(M/K^\Delta)\,
|\ \frak{L}|\frak{l}\,\rangle.
\end{equation}

\

Denote by $U_{F,\frak{L}}(p)$ be the pro-$p$-part of the local unit group at the prime $\frak{L}$
of a number field $F$ of finite degree.
Define
\begin{equation}\label{U}
\mathcal{U}_{K^{\Delta},\frak{l}}(p):=\varprojlim_{k\subseteq F\subseteq K^{\Delta}, [F:\Q]<\infty}
\prod_{\frak{L}\in S_{\frak{l}}(F)}U_{F,\frak{L}}(p)
\end{equation}
to be the projective limit of the pro-$p$-part of the semi-local unit groups of $F$ at $\frak{l}$
with respect to the norm maps, 
where $S_\frak{l}(F)$ stands for the set of all the primes of $F$ lying over $\frak{l}$
and $F$ runs over
all the subfields of $K^\Delta$ with $k\subseteq F$ and $[F:\Q]<\infty$.

Denote by $L^\mathrm{ab}_{p,\{\frak{l}\}}(K^{\Delta})/K^\Delta$
the maximal abelian $p$-extension unramified outside $\frak{l}$.
Then class field theory gives the exact sequence
\[
\varprojlim_{F\subseteq K^\Delta,\,[F:\Q]<\infty}E_F\longrightarrow
\mathcal{U}_{K^\Delta,\frak{l}}(p)\longrightarrow
\Gal(L^{\mathrm{ab}}_{p,\{\frak{l}\}}(K^\Delta)/L_p^{\mathrm{ab}}(K^\Delta))
\longrightarrow 1.
\]
Here it follows from Lemma \ref{lem-un}
that $\varprojlim_{F\subseteq K^\Delta,\,[F:\Q]<\infty}E_F=0$.
Hence we get the isomorphism
\begin{equation}\label{cft}
\mathcal{U}_{K^{\Delta},\frak{l}}(p)\simeq
\Gal(L^\mathrm{ab}_{p,\{\frak{l}\}}(K^\Delta)/L^\mathrm{ab}_p(K^\Delta)).
\end{equation}
On the other hand, we see that
\[
L^{\mathrm{ab}}_p(K^\Delta)\subseteq M
\subseteq L^{\mathrm{ab}}_{p,\{\frak{l}\}}(K^\Delta),
\]
and 
$M/L^{\mathrm{ab}}_p(K^\Delta)$
is the maximal subextenion of $L^{\mathrm{ab}}_{p,\{\frak{l}\}}(K^\Delta)/L^{\mathrm{ab}}_p(K^\Delta)$ such that every ramified prime in $M/L_p(K^\Delta)$
has ramification index $p$.
Hence we see that
\begin{equation}\label{YU}
\mathrm{Tor}(Y_\Delta^{(p)})=\langle I_\frak{L}(M/K^\Delta)\,
|\ \frak{L}|\frak{l}\,\rangle
=\Gal(M/L^{\mathrm{ab}}_p(K^\Delta))
\simeq
\mathcal{U}_{K^{\Delta},\frak{l}}(p)/p
\end{equation}
by using \eqref{tor} and \eqref{cft}.

$U_{F,\frak{L}}(p)/p$ is a cyclic group of oder $p$ on which
$D_\frak{l}(K/k)$ acts on trivially
since $N(\frak{l})\equiv 1\pmod{p}$.
Hence it follows from \eqref{U} that
\[
\mathcal{U}_{K^{\Delta},\frak{l}}(p)/p
\simeq\F_p[[\Gal(\Gal(K^\Delta/k)/
D_\frak{l}(K^\Delta/k))]]
\simeq\F_p[[\Gal(K/k)/D_\frak{l}(K/k))]]
\]
as $\Gal(K/k)$-modules, noting that
$\Delta\subseteq D_\frak{l}(K/k)$.
Thus the assetion of the lemma follows
from \eqref{YU}.
\hfill$\Box$

\

{\bf Proof of Proposition 2}\ \ \ 
Let $N/\Q$ be the Galois closure
of $k_1k_2/\Q$.
We choose an odd prime number $p$ such that $\Q(\mu_p)\cap N=\Q$.
For any given prime number $q\ne p$ unramified in $N$, there exists a prime number $l$
such that $[l,N/\Q]=[q,N/\Q]$, $l\equiv 1\pmod {p}$,
and $l$ is unramified in $k_1k_2$ by the \v{C}ebotarev density theorem, where $[r,N/\Q]$ stands for the Frobenius conjugacy class of $r$ in $\Gal(N/\Q)$ for any 
prime number $r$.
Then we have
\begin{equation*}\label{g_i}
g_{k_1}(q)=g_{k_1}(l),\ g_{k_2}(q)=g_{k_2}(l).
\end{equation*}
Hence if $g_{k_1}(l)=g_{k_2}(l)$ holds for all the 
prime numbers $l\equiv 1\pmod{p}$ unramified in $k_1k_2/\Q$, then 
$g_{k_1}(q)=g_{k_2}(q)$ for all but finitely many prime numbers $q$,
which in turn implies
$k_1\approx k_2$
by Theorem C.

We write $k$ and $K$
for $k_i$ and $K_i$, respectively.
Let $l\equiv 1\pmod{p}$ be a prime number
unramified in $k/\Q$ and $\frak{l}$
a prime of $k$ lying over $l$.
Then it follows from Lamma \ref{lem-DI}
that $I_\frak{l}(K/k)$ has the subgroup
$\Delta_\frak{l}$ of order $p$.
Then it follows from Lemma \ref{lem-Y}
that
\begin{equation}\label{stab}
\mathrm{Stab}_{\Gal(K/k)}
(\mathrm{Tor}(Y_{{\Delta}_\frak{l}}^{(p)}))
:=\{\alpha\in\Gal(K/k)\,|\,\alpha y=y
\mbox{ for all $y\in Y_{\Delta_\frak{l}}$}
\}
=D_{\frak{l}}(K/k).
\end{equation}
We note that
the $\Gal(K/k)$-module structure of $Y_\Delta^{(p)}$
is determined only by   
the group structure of $\Gal(L^{\ab}(K)/k)$
and the subgroup
$\Delta\subseteq\Gal(K/k)$.
The prime number $l$ is charactelized from the structure of $D_\frak{l}(K/k)$
by the unique prime number $l$ such that
$\mathrm{rank}_{\Z_l}D_\frak{l}(K/k)=2$
by Lemma \ref{lem-DI}.
Hence we find from Lemmas \ref{lem-sep}, \ref{lem-Y}
and \eqref{stab} that
\begin{align*}
&g_k(l)\\
&\!=\!\#\{\Delta\subseteq\Gal(K/k)\,|\,
\Delta\simeq\Z/p,\,\mathrm{Tor}(Y_\Delta^{(p)})\ne 0,\, \mathrm{rank}_{\Z_l}
\mathrm{Stab}_{\Gal(K/k)}(\mathrm{Tor}(Y_\Delta^{(p)}))=2\}
\end{align*}
for every prime number $l\equiv 1\pmod{p}$ which is unramified in $k/\Q$.
This means $g_k(l)$ is determined by the group structure of $\Gal(L^\mathrm{ab}(K)/k)$.

Therefore, it follows from the isomorphism
$
\Gal(L^{\mathrm{ab}}(K_1)/k_1)
\simeq\Gal(L^{\mathrm{ab}}(K_2)/k_2)
$
that $g_{k_1}(l)=g_{k_2}(l)$ holds
for every prime number $l\equiv 1\pmod{p}$,
unramified in $k_1k_2/\Q$,
from which we deduce 
$g_{k_1}(q)=g_{k_2}(q)$ for 
all but finitely many prime numbers $q$ as we have seen.

Now, by using Theorem C, we conclude that
$k_1\approx k_2$.
Thus we have proved Proposition 2, from which Proposition 1 follows.
\hfill$\Box$

\

\section{Construction of field isomorphisms}
Let $k_i$ be a number field of finite degree, and put $\tilde{k_i}=k_i(\mu_\infty)$ for $i=1,2$.
Denote by $L_i$ the maximal unramified extension field of $\tilde{k}_i$.

In this section, assuming the existence of
an isomorphism 
$\varphi:\Gal(L_1/k_1)\simeq\Gal(L_2/k_2)$
of pro-finite groups,
we construct a field isomorphism
\[
\tau:L_1\simeq L_2
\]
such that $\tau(k_1)=k_2$ and
$\varphi(x)=\tau x\tau^{-1}$
for every $x\in\Gal(L_1/k_1)$.
We are based on Proposition 1 given in the precedent section and the method of Uchida\cite{Uch79}
to construct a field isomorphism $\tau$.

We first recall the following facts
on arithmetically equivalence:
\begin{thme}(cf.\cite[Theorems (1.3),(1.4),(1.6)]{Kli})\ \ 
Let $F_1$ and $F_2$ be number fields of finite degree such that $F_1\approx F_2$.
Then the followings holds:

(1)\ \ \ $F_1$ and $F_2$ has the common Galois closure over $\Q$.

(2)\ \ \ For any finite Galois extension $F_0/\Q$, we have $F_0F_1\approx F_0F_2$.

(3)\ \ \ Let $N/\Q$ be a finite Galois extension with $F_1F_2\subseteq N$.
Then there exists a bijection
\[
\gamma:\Gal(N/F_1)\longrightarrow\Gal(N/F_2)
\]
such that 
$\gamma(x)=\tau_x x\tau_x^{-1}$
for each $x\in\Gal(N/F_1)$ with some $\tau_x\in\Gal(N/\Q)$.
\end{thme}

\

Let $K_1/k_1$ be any finite Galois subextension of $L_1/k_1$, and
$K_2$ is the corresponding intermediate
field of $L_2/k_2$ by $\varphi$, namely,
$K_2=L_2^{\varphi(\Gal(L_1/K_1))}$.
Then $K_2/k_2$ is also a Galois extension
and $\varphi$ induces the isomorphism
\[
\varphi_{K_1}:\Gal(K_1/k_1)\simeq\Gal(K_2/k_2).
\]
Let $K/\Q$ be a finite Galois extension
containing $K_1K_2$ and put $G:=\Gal(K/k)$
\begin{lem}\label{lem-free}
Let $p$ be a prime number with 
$p\nmid\# G$ and $r$ a positive integer.
Then there exists a Galois extension
$M/\Q$ containing $K$ such that
$\Gal(M/K)\simeq\F_p[G]^{\oplus r}$
as $G$-modules when we define the
$G$-action on $\Gal(M/K)$
via inner automorphisms of $\Gal(M/\Q)$,
and that the maximal abelian $p$-subextension $M_1/K_1$ of $M/K_1$
is a subextension of $L_1/K_1$.
\end{lem}
{\bf Proof.}\ \ \ 
By the \v{C}ebotarev density theorem,
there exist degree one principal prime ideals $(\Lambda_i)\ (1\le i\le r)$ of $K(\mu_p)$ such that
$\Lambda_i\equiv 1\pmod{p^2}$ and totally positive, and that
$(\Lambda_i)$'s are unramified in $K/\Q$ and lying over distinct rational primes
$l_i$'s with $l_i\equiv 1\pmod{p}$.
Define $M/K$ to be the maximal abelian $p$-subextension
of $M':=K(\mu_p,\,\sqrt[p]{\sigma\Lambda_i}\,|\,\sigma\in\Gal(K(\mu_p)/\Q)
,\ 1\le i\le r)/K$.
Then we have
\begin{align*}
\Gal(M'/K(\mu_p))
&\simeq
\mathrm{Hom}_{\F_p}(\F_p[\Gal(K(\mu_p)/\Q)]^{\oplus r},\mu_p)\\
&\simeq\mathrm{Hom}_{\F_p}(\F_p[\Gal(K(\mu_p)/\Q)](-1),\F_p)^{\oplus r}
\end{align*}
as $\Gal(K(\mu_p)/\Q)$-modules by Kummer
duality, $(-1)$ denoting the Tate twist.
Hence we see that
\begin{align*}
\Gal(M/K)&\simeq\Gal(M(\mu_p)/K(\mu_p))\simeq
\Gal(M'/K(\mu_p))_{\Gal(K(\mu_p)/K)}\\
&\simeq
\mathrm{Hom}_{\F_p}((\F_p[\Gal(K(\mu_p)/\Q)](-1))^{\Gal(K(\mu_p)/K)},\F_p)
^{\oplus r}.
\end{align*}
Here, 
\[
(\F_p[\Gal(K(\mu_p)/\Q)](-1))^{\Gal(K(\mu_p)/K)}=
\F_p[\Gal(K(\mu_p)/\Q)]\varepsilon
\simeq\F_p[G]
\]
as $G$-modules
holds for
\[
\varepsilon:=\sum
_{\delta\in\Gal(K(\mu_p)/K)}\chi(\delta)\delta^{-1}
\in\F_p[\Gal(K(\mu_p)/\Q)](-1),
\]
where $\chi:\Gal(K(\mu_p)/K)\longrightarrow\F_p^\times$ stands for the cyclotomic character.

Then we see that
$\Gal(M/K)\simeq\F_p[G]^{\oplus r}$
as $G$-modules, and that
the ramified primes in $M/K$ are exactly the primes lying over $l_i$'s whose ramification indexes equal $p$.
Hence, if we denote by $M_1/K_1$ the maximal abelian $p$-subextension of $M/K_1$,
then the ramified primes in $M_1/K_1$  
are lying over $l_i$'s, whose ramification indexes are $p$ because $p\nmid[K:K_1]$.
Since the prime $l_i$ is unramified in $k_1/\Q$ and $l_i\equiv 1\pmod{p}$,
the primes lying over $l_i$'s
are unramified in $M_1(\mu_{l_i})/K_1(\mu_{l_i})$.
Therefore we see that
$M_1(\mu_\infty)/
K_1(\mu_\infty)$ is unramified abelian $p$-extension.
Because $K_1(\mu_\infty)\subseteq L_1$,
we conclude that $M_1\subseteq L_1$.
\hfill$\Box$

\

Let $H_i:=\Gal(K_i/k_i)$ ($i=1,2$) and
$M/K$ an extension given by Lemma \ref{lem-free}
such that
\begin{equation}\label{A}
A:=\Gal(M/K)\simeq\bigoplus_{h\in H_1}
\F_p[G]u_h,
\end{equation}
as $G$-modules, where the right hand side
is the free $\F_p[G]$-modules 
with basis $\{u_h\,|\,h\in H_1\}$.

For each $h\in H_1$,
let $M_h/K$ be the subextension
of $M/K$ such that 
\[
\Gal(M/M_h)
\simeq\bigoplus_{h\ne h'\in H_1}\F_p[G]u_{h'}
\]
under isomorphism \eqref{A}.
We note that $M_h/\Q$ is a Galois extension.

Define $M_{1,h}/K_1$ to be the maximal abelian $p$-subextension of $M_h/K_1$ which is abelian over
$K_1^{\langle h\rangle}$ for each $h\in H_1$. 
By our choice of $M$, we see that $M_{1,h}\subseteq L_1$.
Let $M_{2,h}$ be 
the field coresponding to
$M_{1,h}$ by $\varphi_{K_1}$ for each $h\in H_1$.
Then it follows from Proposition 1
that $M_{1,h}\approx M_{2,h}$.
This means $M_{2,h}$ is a subfield of
the Galois closure of $M_{1,h}$ over $\Q$
by Theorem E (1),
which is a subfield of $M_h$ since $M_h/\Q$ is Galois.
Hence we see that $K_2\subseteq M_{2,h}\subseteq M_h$  
and that $\varphi$ induces
$\Gal(M_{2,h}/K_2)\simeq
\Gal(M_{1,h}/K_1)$,
which are abelian $p$-groups, and
$\Gal(M_{2,h}/K_2^{\langle\varphi_{K_1}(h)\rangle})\simeq
\Gal(M_{1,h}/K_1^{\langle h\rangle})$,
which is abelian.

Furthemere, for any abelian $p$-subextension $F_2/K_2$ of $M_h/K_2$
which is abelian over $K_2^{\langle\varphi_{K_1}(h)\rangle}$,
we see that the corresponding field $F_1$ to $F_2$ by $\varphi^{-1}$ is an intermediate field of $M_h/K_1$
by a similar argument above, which is $p$-abelian over $K_1$, and
$\Gal(F_1/K_1^{\langle h\rangle})
\simeq\Gal(F_2/K_2^{\langle\varphi_{K_1}(h)\rangle})$ is abelian.
Hence we have $F_1\subseteq M_{1,h}$, which implies $F_2\subseteq M_{2,h}$.
Therefore we conclude that 
$M_{2,h}/K_2$ is the maximal abelian $p$-subextension of $M_h/K_2$ which is abelian over $K_2^{\langle\varphi_{K_1}(h)\rangle}$.

\

\begin{lem}\label{lem-A}
(1)\ \ \ We have
\[
K\prod_{h\in H_1}M_{1,h}\approx K\prod_{h\in H_1}M_{2,h}.
\]

(2)\ \ \ Let $N_i:=\Gal(K/K_i)\ (i=1,2)$.
Then we have
\begin{align*}
A_1:=\Gal(M/K\prod_{h\in H_1}M_{1,h})
\simeq\bigoplus_{h\in H_1}
(I_{N_1} +(\overline{h}-1)\F_p[G])u_h,\\
A_2:=\Gal(M/K\prod_{h\in H_1}M_{2,h})
\simeq\bigoplus_{h\in H_1}
(I_{N_2} +(\overline{\varphi_{K_1}(h)}-1)\F_p[G])u_h,
\end{align*}
where $I_{N_i}:=\sum_{n\in N_i}(n-1)\F_p[G]$, 
$\overline{h}\in\Gal(K/k_1)$ and  
$\overline{\varphi_{K_1}(h)}\in\Gal(K/k_2)$
are lifts of $h$ and $\varphi_{K_1}(h)$, respectively.
\end{lem}
{\bf Proof.}\ \ \ 
(1)\ \ \ Because $M_{2,h}$ is corresponding to $M_{1,h}$ by $\varphi$ for each $h\in H_1$,
the coposite field 
$\prod_{h\in H_1}M_{2,h}$
is corresponding to $\prod_{h\in H_1}M_{1,h}$ by $\varphi$.
Then it follows from Proposition 1 that $\prod_{h\in H_1}M_{1,h}\approx\prod_{h\in H_1}M_{2,h}$. Hence the assertion follows from
Theorem E(2) since $K/\Q$ is Galois.

\

(2)\ \ \ 
Recall that $M_{1,h}/K_1$ and
$M_{2,h}/K_2$ are the maximal
abelian $p$-subextensions of $M_h/K_1$
and $M_h/K_2$ such that $h$ and
$\varphi_{K_1}(h)$ act trivially on
$\Gal(M_{1,h}/K_1)$ and 
$\Gal(M_{2,h}/K_2)$, respectively.
Hence we find that 
$N_1$ and $\overline{h}$
acts trivially on $
\Gal(M/K)/\Gal(M/K{M_{1,h}})
\simeq\Gal(KM_{1,h}/K)
\simeq\Gal(M_{1,h}/K_1)$,
which implies
\[
J_h:=\bigoplus_{h\ne h'\in H_1}\F_p[G]u_{h'}
\oplus(I_{N_1}+(\overline{h}-1)\F_p[G])u_h\subseteq
\Gal(M/KM_{1,h})
\]
It follows from the definition of $J_h$ that $\Gal(M^{J_h}/K_1)$ has a direct factor
naturally isomorphic to $N_1=\Gal(K/K_1)$ since $p\nmid\# N_1$.
Hence there is an abelian $p$-subextension
$F/K_1$ such that $KF=M^{j_h}$ and $F\cap K=K_1$,
and we see that $F/K_1^{\langle h\rangle}$
is abelian by the definition of $J_h$. 
This means $F\subseteq M_{1,h}$ and
$\Gal(M/KM_{1,h})\subseteq J_h$.
Thus we conclude that $J_h=\Gal(M/KM_{1,h})$
and 
\[\Gal(M/K\prod_{h\in H_1}M_{1,h})
=\bigcap_{h\in H_1}J_h=\bigoplus_{h\in H_1}
(I_{N_1} +(\overline{h}-1)\F_p[G])u_h.
\]
We obtain the assertion also for $A_2$ by the similar way.
\hfill$\Box$

\

\begin{lem}\label{lem-tau}
There exists $\tau_0\in G$ and $m_h\in\Z$
for each $h\in H_1$ such that
$\tau_0 N_1\tau_0^{-1}=N_2$, and
\[
(\tau_0\overline{h}\tau_0^{-1})|_{K_2}=
\varphi_{K_1}(h)^{m_h}\in
\Gal(K_2/k_2)
\]
holds for every $h\in H_1$, 
where $\overline{h}\in\Gal(K/k_1)$ is a lift of $h$. 
\end{lem}
{\bf Proof.}\ \ \ 
By using Lemma \ref{lem-A} (2), we define
\[
\alpha:=\sum_{n\in N_1}(n-1)u_1+
\sum_{h\in H_1-\{1\}}(\overline{h}-1)u_h
\in A_1.
\]
Then it follows from Lemma \ref{lem-A} (1) and Theorem E (3) that there exists $\tau_0\in G$
such that
\begin{equation}\label{taualpha}
\tau_0\cdot\alpha\in A_2,
\end{equation}
where `` $\cdot$ " denotes the $G$-action on $A$.
We derive from \eqref{taualpha} and Lemma \ref{lem-A} (2) that
\begin{equation}\label{t1}
\tau_0\sum_{n\in N_1}(n-1)\in I_{N_2},
\end{equation}
and
\begin{equation}\label{th}
\tau_0(\overline{h}-1)\in I_{N_2}
+(\overline{\varphi_{K_1}(h)}-1)\F_p[G],
\end{equation}
for $h\in H_1-\{1\}$.

By operating 
$T_{N_2}:=\sum_{n\in N_2}n\in\F_p[G]$
on \eqref{t1}, we obtain the equality
\[
T_{N_2}\tau_0\sum_{n\in N_1}n=(\#N_1)T_{N_2}\tau_0,
\]
from which we see that
for each $n_1\in N_1$, there exists $n_2\in N_2$
such that $\tau_0n_1=n_2\tau_0$.
This implies $\tau_0N_1\tau_0^{-1}\subseteq N_2$.
Since $K_1\approx K_2$ by Proposition 1,
we find that $[K_1:\Q]=[K_2:\Q]$ by Theorem E (3),
which implies $\# N_1=\# N_2$.
Thus we conclude that $\tau_0N_1\tau_0^{-1}=N_2$.

By operating $T:=\sum_{t\in\langle N_2,\overline{\varphi_{K_1}(h)}\rangle}t\in
\F_p[G]$ on \eqref{th}, we get
\[
T\tau_0\overline{h}=T\tau_0,
\]
from which we see that
there exists $t\in\langle N_2,\overline{\varphi_{K_1}(h)}\rangle$
such that 
$\tau_0\overline{h}=t\tau_0$.
Then we conclude that
\[
(\tau_0\overline{h}\tau_0^{-1})N_2
=\overline{\varphi_{K_1}(h)}^{m_h}N_2,
\]
which means 
\[
\tau_0\overline{h}\tau_0^{-1}|_{K_2}
=\varphi_{K_1}(h)^{m_h},
\]
for a certain $m_h\in\Z$ for each $h\in H_1$.\hfill$\Box$

\

We note that  
$\Gal(K/\tau_0(k_1))
=\tau_0\Gal(K/k_1)\tau_0^{-1}\subseteq\Gal(K/k_2)$ holds by the above Lemma,
which in turn implies 
\begin{equation}\label{kisom}
\tau_0(k_1)=k_2,\ \tau_0\Gal(K/k_1)\tau_0^{-1}=\Gal(K/k_2),
\end{equation}
because $[K:k_1]=[K:k_2]$ holds by $k_1\approx k_2$
and Theorem E. Hence $k_1\simeq k_2$ holds.

\

In what follows we will show that the assertion of Lemma \ref{lem-tau} holds for $m_h=1$ in fact. 
\begin{lem}\label{lem-GalM}
Let $K_1/k_1$ be any finite Galois subextension of
$L_1/k_1$ and $p$ a prime number. 
Then there exists a Galois subextension $M_1/k_1$
of $L_1/k_1$ with $K_1\subseteq M_1$
such that
\[
\Gal(M_1/K_1)\simeq\F_p[H_1]
\]
as $H_1=\Gal(K_1/k_1)$-modules, where
the $H_1$-acsion on $\Gal(M_1/K_1)$ is defined via
inner automorphisms of $\Gal(M_1/k_1)$.
\end{lem}
{\bf Proof.}\ \ \ We can show the existence of $M_1$ in a similar way to the proof of Lemma \ref{lem-free}. 
Choose a principal degree one prime ideal $(\Lambda)$ of $K_1(\mu_p)$ such that
$\Lambda\equiv 1\pmod{p^2}$ and totally positive, and that
the rational prime $l$ below $(\Lambda)$ is unramified in $k_1/\Q$.
Then the maximal abelian $p$-subextension
$
M_1/K_1$ of 
$K_1(\mu_p)(\sqrt[p]{\sigma\Lambda}\,|\,\sigma\in\Gal(K_1(\mu_p)/k_1))/K_1$ satisfies our requirement.
\hfill$\Box$

\

Now we will give the following crucial
proposition:
\begin{prop}
Let $K_1/k_1$ be a finite Galois subextension of $L_1/k_1$, and
$K_2\subseteq L_2$ the corresponding field
to $K_1$ by $\varphi$.
Then there exists a field isomorphism
$\tau_{K_1}:K_1\simeq K_2$ such that
$\tau_{K_1}(k_1)=k_2$ and 
\[
\varphi_{K_1}(x)=\tau_{K_1}x\tau_{K_1}^{-1}
\]
for every $x\in\Gal(K_1/k_1)$.
\end{prop}
{\bf Proof.}\ \ \ 
$M_1/K_1$
be a subextension of $L_1/K_1$ given by 
Lemma \ref{lem-GalM}.
Let $K_2$ and $M_2$ be the fields corresponding
to $K_1$ and $M_1$ by $\varphi$, respectively.
Then $\varphi$ induces the isomorphism
\begin{equation}\label{phiM}
\varphi_{M_1}:\Gal(M_1/k_1)\simeq\Gal(M_2/k_2)
\end{equation}
with $\varphi_{M_2}(\Gal(M_1/K_1))=\Gal(M_2/K_2)$.
Let $M/\Q$ be a finite Galois extension
containing $M_1M_2$.
Then by applying Lemma \ref{lem-tau} to $M_1,\ M_2$,
and $M$
as $K_1,\ K_2$, and $K$, respectively,
we get $\tau_0\in\Gal(M/\Q)$ such that
\begin{equation}\label{GM}
\tau_0\Gal(M/M_1)\tau_0^{-1}=\Gal(M/M_2) \end{equation}
and 
\begin{equation}\label{conj}
(\tau_0\overline{x}\tau_0^{-1})|_{M_2}
=\varphi_{M_1}(x)^{m_x}
\end{equation}
holds for each $x\in\Gal(M_1/k_1)$ with a certain
$m_x\in\Z$,
where $\overline{x}\in\Gal(M/k_1)$
is a lift of $x$. 
Here we note that 
$\tau_0\Gal(M/k_1)\tau_0^{-1}=\Gal(M/k_2)$
holds by \eqref{kisom}.
Put $A:=\Gal(M_1/K_1)=\F_p[\Gal(K_1/k_1)]u$,
$u$ being a basis of the free $\F_p[\Gal(K_1/k_1)]$-module $A$.
Then we have 
\[
\Gal(M_2/K_2)=\varphi_{M_1}(\Gal(M_1/K_1))=
\F_p[\Gal(K_2/k_2)]\varphi_{M_1}(u)\simeq\F_p[\Gal(K_2/k_2)],
\]
as $\Gal(K_2/k_2)$-modules
by \eqref{phiM}.

Let $x\in\Gal(M_1/k_1)$ be any
element.
Then we have
\begin{equation}\label{phiM1}
\begin{aligned}
&\varphi_{M_1}(xux^{-1})^{m_u}\\
&=\varphi_{M_1}(x)
\varphi_{M_1}(u)^{m_u}
\varphi_{M_1}(x)^{-1}\\
&=\varphi_{M_1}(x)\left((\tau_0\overline{u}\tau_0^{-1})|_{M_2}\right)
\varphi_{M_1}(x)^{-1},
\end{aligned}
\end{equation}
where $u\in A$ is the free basis,
by \eqref{conj}.

We also get
\begin{equation}\label{phiM2}
\begin{aligned}
&\varphi_{M_1}(xux^{-1})^{m_{xux^{-1}}}
\\
&=\left(\tau_0(\overline{xux^{-1}})\tau_0^{-1}\right)|_{M_2}\\
&=(\tau_0\overline{x}\tau_0^{-1})|_{M_2}
(\tau_0\overline{u}\tau_0^{-1})|_{M_2}(\tau_0x\tau_0^{-1})^{-1}|_{M_2}
\end{aligned}
\end{equation}

Because we have $\tau_0\overline{u}\tau_0^{-1}|_{M_2}\ne 1$ by \eqref{GM} and
$\varphi_{M_1}(xux^{-1})\in\Gal(M_2/K_2)\simeq\F_p[\Gal(K_2/k_2)]$,
we see that $m_u$ and $m_{xux^-1}$
is prime to $p$ by \eqref{phiM1}
and \eqref{phiM2}.
Furthermore, since we have
$
\varphi_{M_1}(u)^{m_u}=(\tau_0\overline{u}\tau_0^{-1})|_{M_2}
$
from \eqref{conj}, and $p\nmid m_u$.
we see that
\[
(\tau_0\overline{u}\tau_0^{-1})|_{M_2}
\in\Gal(M_2/K_2)=\F_p[\Gal(K_2/k_2)]
(\varphi_{M_1}(u))
\]
is a free $\F_p[\Gal(K_2/k_2)]$-basis of
$\Gal(M_2/K_2)$.

We derive from \eqref{phiM1} and \eqref{phiM2} that
\[
m_{xux^{-1}}\varphi_{M_1}(x)|_{K_2}\cdot(\tau_0\overline{u}\tau_0^{-1})
|_{M_2}
=
m_u(\tau_0\overline{x}\tau_0^{-1})|_{K_2}
\cdot
(\tau_0\overline{u}\tau_0^{-1})_{M_2}
\in\Gal(M_2/K_2),
\]
where `` $\cdot$ " stands for the $\F_p[\Gal(K_2/k_2)]$-action on $\Gal(M_2/K_2)$.

Because  $(\tau_0\overline{u}\tau_0^{-1})_{M_2}
\in\Gal(M_2/K_2)\simeq\F_p[\Gal(K_2/k_2)]$
is a free $\F_p[\Gal(K_2/k_2)]$-basis
of $\Gal(M_2/K_2)$ and $p\nmid m_um_{xux^{-1}}$, we conclude that
\begin{equation}\label{phiK}
\varphi_{K_1}(x|_{K_1})=\varphi_{M_1}(x)|_{K_2}=(\tau_0\overline{x}\tau_0^{-1})
|_{K_2}\in\Gal(K_2/k_2)
\end{equation}
for every $x\in\Gal(M_1/k_1)$.

It follows from \eqref{phiK} that
$
\tau_0\Gal(M/K_1)\tau_0^{-1}\subseteq\Gal(M/K_2)
$, from which
we have 
\[
\tau_0\Gal(M/K_1)\tau_0^{-1}=\Gal(M/K_2),
\]
since $\#\Gal(M/K_1)=\#\Gal(M/K_2)$
by Proposition 1 and Theorem E (3).
This means $\tau_0(K_1)=K_2$.
Furthermore, $\tau_0(k_1)=k_2$ also holds by \eqref{kisom}. 

Now we define the field isomorphism
$\tau_{K_1}:K_1\simeq K_2$ by
$\tau_{K_1}(\alpha):=\tau_0(\alpha)$
for $\alpha\in K_1$.
Then $\tau_{K_1}(k_1)=k_2$
and
\[
\varphi_{K_1}(x)=\tau_{K_1}x\tau_{K_1}^{-1}
\]
for every $x\in\Gal(K_1/k_1)$ holds by \eqref{phiK}.
This completes the proof of the proposition.\hfill$\Box$

\

\section{Proof of Theorem 1}
Now we will giva a proof of Theorem 1.
Let $k_i$ be number field of finite degree
and $L_i$ the maximal unramified extension
of $\tilde{k_i}=k_i(\mu_\infty)$ ($i=1,2$).
Assume that 
\[
\varphi:\Gal(L_1/k_1)\simeq\Gal(L_2/k_2)
\]
is an isomorphism of pro-finite groups.
We will show that there exists the unique
field isomorphism
\[
\tau:L_1\simeq L_2
\]
such that 
\begin{equation}\label{req}
\tau(k_1)=k_2,\ \ 
\varphi(x)=\tau x\tau^{-1}
\end{equation}
holds for every $x\in\Gal(L_1/k_1)$.
This implies the assertion of Theorem 1.

For a finite Galois subextension $K_1/k_1$ of $L_1/k_1$,
let $K_2$ be the intermediate field of $L_2/k_2$ corresponding to $K_1$
by $\varphi$, and define
$T_{K_1}$ to be the set of all the field isomorphisms
\[
\tau_{K_1}:K_1\simeq K_2
\]
such that
\begin{equation}\label{reqK_1}
\tau_{K_1}(k_1)=k_2,\ \ \varphi_{K_1}(x)=\tau_{K_1}x\tau_{K_1}^{-1}
\end{equation}
holds for every $x\in\Gal(K_1/k_1)$,
where $\varphi_{K_1}:\Gal(K_1/k_1)\simeq\Gal(K_2/k_1)$ is the isomorphism inducedc by $\varphi$.

Then it follows from Proposition 3
that $T_{K_1}\ne\emptyset$ for each $K_1$.
Furhermore, if $k_1\subseteq K_1\subseteq K'_1\subseteq L_1$ are inermediate fields finite Galois over $k_1$,
we obtain the map
\[
T_{K'_1}\longrightarrow T_{K_1}
\]
by $\tau_{K'_1}\mapsto\tau_{K'_1}|_{K_1}$.
Since $T_{K_1}$ is a non-empty finite set,
we see that the projective limit
$T$ of $T_{K_1}$'s with respect to
the above maps is not empty, where
$K_1$ runs over all the intermediate fields of $L_1/k_1$ such that $K_1/k_1$ is finite Galois.
Take $(\tau_{K_1})_{K_1}\in T$ and
define the map
\[
\tau:L_1\longrightarrow L_2
\]
by $\tau(\alpha):=\tau_{K_1}(\alpha)$
for each $\alpha\in L_1$ and finite Galois subextension $K_1/k_1$ of $L_1/k_1$
with $\alpha\in K_1$.
We see that $\tau$ is a well-defined
field isomorphism $L_1\simeq L_2$
and satisfies our requirement \eqref{req} by \eqref{reqK_1}.
Thus we have proved the existence of $\tau:L_1\simeq L_2$ with \eqref{req}.

Finally we will show the uniquness of $\tau$ with \eqref{req} in what follows.
We need the following:
\begin{lem}\label{lem-Z}
Let $k_0\subseteq k_1$ be a subfield
such that $L_1/k_0$ is a Galois extension.
Then the centralizer 
$Z_{\Gal(L_1/k_0)}(\Gal(L_1/k_1))$
of $\Gal(L_1/k_1)$ in $\Gal(L_1/k_0)$
is trivial.
\end{lem}
{\bf Proof.}\ \ \ 
Assume that $z\in Z_{\Gal(L_1/k_0)}(\Gal(L_1/k_1))$.
Let $K_1/k_0$ be a finite Galois subextension of $L_1/k_0$
By a similar way to the proof of Lemma \ref{lem-GalM}, we see that there exists finite Galois subextension $M_1/k_0$ of $L_1/k_0$ containing $K_1$
such that $\Gal(M_1/K_1)\simeq
\F_p[\Gal(K_1/k_0)]$
as $\Gal(K_1/k_0)$-modules,
where 
the $\Gal(K_1/k_0)$-action on 
$\Gal(M_1/K_1)$ is given via inner
automorphisms of $\Gal(M_1/k_0)$.
Then $z|_{K_1}\in\Gal(K_1/k_0)$ acts on $\Gal(M_1/K_1)$
trivially, which implies
$z|_{K_1}=1$ since the $\Gal(K_1/k_0)$-action on $\Gal(M_1/K_1)$ is faithful.
Because $K_1$ can be arbitral
finite Galois subextension of $L_1/k_0$, we conclude that $z=1$.
\hfill$\Box$

\

Assume taht 
$\tau_1, \tau_2:L_1\simeq L_2$
are field isomorphisms satisfying 
\eqref{req} for $\tau=\tau_1,\tau_2$.
Put
$z:=\tau_1^{-1}\tau_2\in
\Gal(L_1/k_0)$, where $k_0$ be the minimal subfield of $L_1$ such that
$L_1/k_0$ is Galois.
Then  we have
\[
zxz^{-1}=\tau_1^{-1}(\tau_2 x\tau_2^{-1})\tau_1
=\tau_1^{-1}\varphi(x)\tau_1
=\varphi^{-1}\varphi(x)=x
\]
for any $x\in\Gal(L_1/k_1)$
by our assumption.
Hence we conclude that 
\[
z\in Z_{\Gal(L_1/k_0)}(\Gal(L_1/k_1))=1
\]
by Lemma \ref{lem-Z},
which implies $\tau_1=\tau_2$.
Thus we have shown the uniqueness of $\tau$ with \eqref{req}.
This completes the proof of Theorem 1.
\hfill$\Box$
\section{Remarks}
1.\ \ In the case of a curve $C$ over a finite field $\F$,
it is known that
$\Gal(L(C)/\overline{\F}(C))$
is a characteristic subgroup of
$\Gal(L(C)/\F(C))$ (see \cite[Proposition(3.3)]{Tam}).
Furthermore, for every $\varphi\in\mathrm{Aut}(\Gal(L(C)/F(C)))$, the induced automorphism of
$G_{\F}\simeq\Gal(\overline{\F}(C)/\F(C))$ by $\varphi$ is the identity (see \cite[Proposition(3.4)]{Tam}). 

In the case of a number field $k$, analogous assertions also hold, namely,
it follows from Theorem 1 that 
$\Gal(L(\tilde{k})/\tilde{k})$ is a characteristic subgroup of
$\Gal(L(\tilde{k})/k)$
and that for every $\varphi\in\mathrm{Aut}(\Gal(L(\tilde{k})/k))$, the induced automorphism of
$\Gal(\tilde{k}/k)$ by $\varphi$ is the identity.

\

2.\ \ The assertion of Theorem 1 holds
even if we replace $L_i=L(\tilde{k_i})$ with
any solvably closed unraramified extension field $L_i'$ over $\tilde{k_i}$
which is Galois over $k_i$, namely, there is no non-trivial unramified abelian extensions over $L_i'$.
The maximal unramified solvable extension
over $\tilde{k_i}$ is an example of such 
$L_i'$.

Indeed, because the maximal abelian subextenions
of $L_i/k_i$ and $L'_i/k_i$ coincide,
and $L^{\mathrm{ab}}(F)\subseteq L_i'$,
holds for any $k_i\subseteq F\subseteq L_i'$, the assertion of Proposition 1
holds even if we replace $L_i$ with $L_i'$.
Furthemore, for a finite Galois subextension $K_1/k_1$ of $L_1'/k_1$, 
if an abelian extension $M_1/K_1$ satisfies that
$\tilde{M_1}/\tilde{K_1}$ is unramified,
then $M_1\subseteq L_1'$ holds.
Therefore, the arguments of sections 3 and 4 work for $L_i'/k_i$.

\

\

\noindent
Manabu Ozaki,\par\noindent
Department of Mathematics,\par\noindent
School of Fundamental Science and Engineering,\par\noindent
Waseda University,\par\noindent
Ohkubo 3-4-1, Shinjuku-ku, Tokyo, 169-8555, Japan\par\noindent
e-mail:\ \verb+ozaki@waseda.jp+
\end{document}